\documentclass [11pt]{article}
\usepackage{graphicx}
\newtheorem{theorem}{Theorem}
\newtheorem{lemma}[theorem]{Lemma}
\newtheorem{proposition}[theorem]{Proposition}

\newtheorem{corollary}[theorem]{Corollary}

\newtheorem{definition}[theorem]{Definition}

\newtheorem{notation}[theorem]{Notation}
\newtheorem{remark}[theorem]{Remark}
\newtheorem{example}[theorem]{Example}

\begin{document}

\title{{\bf On syzygies of Segre embeddings}}
\author{{\bf  Elena Rubei }}
\date{}
\maketitle

\begin{abstract} 
We study the syzygies  of the ideals of the Segre embeddings.
Let $d \in {\bf N}$, $ d \geq 3$;
we prove that the line bundle ${\cal O}(1,...,1)$
on the $P^1 \times .... \times P^1 $ ($d$ copies) satisfies
Property $N_p$ of Green-Lazarsfeld if and only if $p \leq 3$. 
Besides we prove that if we have a projective variety not satisfying
Property $N_p$ for some $p$, then the product of it with any other projective
variety does not satisfy Property $N_p$. From this we deduce also other 
corollaries about syzygies of Segre embeddings.
\end{abstract}

\def\thefootnote{}
\footnotetext{
{\bf Address:} Elena Rubei, Dipartimento di Matematica ``U. Dini'',
via Morgagni 67/A, 
50134 Firenze, Italia

\hspace{0.24cm}{\bf E-mail address:} rubei@math.unifi.it

\hspace{0.24cm}{\bf 2000 Mathematical Subject Classification:} 
14M25, 13D02.


}

\def\thefootnote{\arabic{footnote}}

\setcounter{footnote}{0}

\section{Introduction}

Let  $ M $ be a very ample  line bundle on a smooth
complex projective variety $Y$ and let  $\varphi_{M}: Y \rightarrow
{\bf P}(H^{0}(Y, M)^{\ast})$ be  the map  associated to $M$.
We  recall the definition of Property $N_{p}$ of Green-Lazarsfeld,
 studied for the first time  by Green in  \cite{Green1}
(see also \cite{G-L}, \cite{Green2}):

{\em
let $Y$ be a smooth complex projective variety and let $L$ be a very ample
 line bundle on $Y$ defining an embedding $\varphi_{L}: Y \hookrightarrow
{\bf P}={\bf P}(H^{0}(Y,L)^{\ast })$;
set $S= S(L)= Sym^{\ast}H^{0}(L) $,
the homogeneous coordinate ring of the projective space ${\bf P}$, and
consider
the graded $S$-module  $G=G(L)= \oplus_{n} H^{0}(Y, L^{n})$; let $E_{\ast}$
\[ 0 \longrightarrow E_{n} \longrightarrow E_{n-1} \longrightarrow 
... \longrightarrow E_{0} \longrightarrow G \longrightarrow 0\] 
be a minimal graded free resolution of $G$;
the line bundle $L$ satisfies Property
$ N_{p}$ ($p \in {\bf N}$) if and only if

\hspace{1cm} $E_{0}= S$

\hspace{1cm} $E_{i}= \oplus S(-i-1)$    \hspace{1cm}
for $1 \leq i \leq p $.}

(Thus $L$ satisfies Property $N_{0}$  if and only if  $Y \subset
{\bf P}(H^{0}(L)^{\ast })$ is projectively normal, i.e.  $L$
is normally generated; $L$ satisfies Property $N_{1}$
if and only if $L$ satisfies Property $N_{0}$
and the homogeneous ideal $I$ of $Y \subset
{\bf P}(H^{0}(L)^{\ast })$ is generated by quadrics;
$L$ satisfies Property $N_{2}$  if and only if
$L$ satisfies Property $N_{1}$ and the module of syzygies among
quadratic
generators $Q_{i} \in I$ is spanned by relations of the form
 $\sum L_{i}Q_{i}=0$, where $L_{i}$  are linear polynomials;
and so on.)

Now let $L= 
{\cal O}_{{\bf P}^{n_1}\times... \times {\bf P}^{n_d}}(a_1, ...,a_d)$, 
where $d,a_1,..., a_d, n_1,...,n_d$ are positive integers. 
The known results on the syzygies in this case are the following:

{\it Case $d=1$}, i.e. the case of the Veronese embedding:

\begin{theorem} \label{Gr} ({\bf Green}) \cite{Green1}.
Let $a$ be a positive integer.
The line bundle ${\cal O}_{{\bf P}^{n}} (a)$ 
  satisfies Property $N_a$.
\end{theorem}

\begin{theorem} \label{OP}  ({\bf Ottaviani-Paoletti}) \cite{O-P}.
If $n \geq 2$,
  $a \geq 3$ and the bundle
${\cal O}_{ {\bf P}^{n} } (a)$ satisfies Property $N_{p} $, then $
p \leq 3a -3$.
\end{theorem}

\begin{theorem} \label{JPW} ({\bf Josefiak-Pragacz-Weyman}) \cite{J-P-W}.
The bundle
${\cal O}_{ {\bf P}^{n} } (2)$ satisfies Property $N_{p} $ if and only if $
p \leq 5$ when $n \geq 3$ and  for all $p$ when $n=2$.
\end{theorem}

(See \cite{O-P} for a more complete bibliography.)

{\it
Case $d=2$:}

\begin{theorem} \label{GP} ({\bf Gallego-Purnapranja}) \cite{G-P}.
Let $a,b \geq 2$. 
The line bundle ${\cal O}_{{\bf P}^{1} \times {\bf P}^{1} }(a,b)$ 
 satisfies Property $N_p$ if and only if $p \leq 2a+2b-3 $.
\end{theorem}

\begin{theorem} \label{LPW} ({\bf Lascoux-Pragacz-Weymann}) \cite{Las}, 
\cite{P-W}. Let $n_1, n_2 \geq 2$.
   The line bundle ${\cal O}_{{\bf P}^{n_1} \times {\bf P}^{n_2}}(1,1)$ 
 satisfies Property $N_p$ if and only if $p \leq 3$.
\end{theorem}

Here we consider 
 $ {\cal O}(1, ...,1)$ on ${\bf P}^{1}\times... \times {\bf P}^{1}$
($d$ times, for any $d$). 
We prove (Section 2):

\begin{theorem} \label{N2}.
The line bundle  $ {\cal O}(1, ...,1)$ on ${\bf P}^{1}\times... \times 
{\bf P}^{1}$
($d$ times) satisfies Property $N_{3}$ for any $d$.
\end{theorem}

Besides we prove (Section 3):

\begin{proposition} \label{prodotti}.
Let $X$ and $Y$ be two projective varieties and let
$L$ be a line bundle on $X$ and $M$ a line bundle on $Y$.
Let $\pi_X : X \times Y \rightarrow X$ and $\pi_Y : X \times Y \rightarrow Y$
be the canonical projections. 
Suppose $L$ and $M$ satisfy Property $N_1$. Let $p \geq 2$.
If $L$ does not satisfy Property $N_p$, then  $\pi_{X}^{\ast}L
\otimes \pi_{Y}^{\ast}M$ does not satisfy Property $N_p$, either. 
\end{proposition}

\begin{corollary} \label{N4}.
Let $a_1 , ..., a_d $ be positive integers
with $a_1 \leq  a_2 \leq ...\leq a_d$.
Suppose $k= max \{i | a_i =1 \}$.
If $k \geq 3$ the line bundle ${\cal O}_{{\bf P}^{n_1} 
\times ...\times {\bf P}^{n_d} }(a_{1},..., a_{d})$ 
does not satisfy Property $N_4$ and  if $d - k \geq 2$
 it does not satisfy Property $N_{2a_{k+1}+2a_{k+2} -2}$.
 \end{corollary}

In particular, from Corollary \ref{N4} and Theorem \ref{N2}, we have:

\begin{corollary} \label{fine}.
Let $ d \geq 3 $. The line bundle ${\cal O}_{{\bf P}^{1} 
\times ...\times {\bf P}^{1} }(1,..., 1)$ ($d$ times)
 satisfies Property $N_p$ if and only if $p \leq 3$.
 \end{corollary}

\section{Proof of Theorem \ref{N2}}

First we have to recall   
some facts on toric ideals from \cite{St}.
       
Let $k \in {\bf N}$. Let $A = \{a_1,..., a_n \}$ be a subset of ${\bf Z}^k$.
The toric ideal  ${\cal I}_A$ is defined as  the  ideal 
in ${\bf C}[x_1, ..., x_n]$
generated as vector space by the
binomials 
$$x_1^{u_1}... x_n^{u_n}
- x_1^{v_1}... x_n^{v_n}$$
 for $(u_1,..., u_n),(v_1,.., v_n) \in {\bf N}^n$, 
with $\sum_{i=1,...,n} u_i a_i = \sum_{i=1,...,n} v_i a_i $. 
           
We have that ${\cal I}_A$ is homogeneous if and only if $\exists \;
\omega \in {\bf Q}^k$ s.t. $\omega  \cdot a_i= 1$ $\forall i = 1,..., n$;
the rings ${\bf C}[x_1, ..., x_n]$ and 
  ${\bf C}[x_1, ..., x_n]/{\cal I}_A$ 
are multigraded by ${\bf N}A$ via $\deg x_i =a_i$;  
the element $x_1^{u_1}...x_n^{u_n}$ 
has multidegree $ b =  \sum_i u_i a_i \in {\bf N}A$ and degree 
$\sum_i  u_i = b \cdot \omega$; we define $\deg b = b \cdot \omega$. 

Theorem 12.12 p.120 in \cite{St}
  studies the syzygies of the
ideal ${\cal I}_A$;  
   for each $b \in {\bf N}A$, let 
$\Delta_b$ be the simplicial complex 
on the set $\{1,..,n \}$  defined   as follows:
$$ \Delta_b = \{F \subset\{1,..,n\}: b - \sum_{i \in F} 
a_i  \in {\bf N}A \}$$
             
(thus, by identifying $\{1,..., n\}$ with $A$, we have: 
$$\Delta_b = \cup_{k \in {\bf N}, a_{i_1},...., a_{i_k} \in A, a_{i_1}+...+
a_{i_k} =b}        < a_{i_1},...., a_{i_k} >,$$
where $< a_{i_1},...., a_{i_k} >$ is the simplex generated by 
$a_{i_1},...., a_{i_k} $):

\begin{theorem} \label{CPS}
 ({\bf Campillo-Pison-Sturmfels})  \cite{C-P},\cite{St} (Thm.
 12.12).
Let $A=\{a_1,...,a_n\}$ 
be a subset of ${\bf Z}^k$ and ${\cal I}_A$ be the associated 
toric ideal.  Let $ 0 \rightarrow E_n \rightarrow ...\rightarrow
 E_1 \rightarrow E_0 \rightarrow G \rightarrow 0$ be a 
minimal free resolution of $G= {\bf C}[x_1,..., x_n]/ {\cal I}_A$
on ${\bf C}[x_1,..., x_n]$.
Each of the generators of $E_j$ has a unique multidegree. 
The number of the generators of 
multidegree $b \in {\bf N}A$ of  $E_{j+1}$  equals the rank of the 
$j$-th
reduced homology group $\tilde{H}_j(\Delta_b, {\bf C}) $ of the 
simplicial 
complex $\Delta_b$.
\end{theorem}

\begin{notation}.
If $\alpha$ is a chain in a topological space, 
$sp(\alpha)$ will denote the support 
of $\alpha $, i.e. the union of the supports of the 
simplexes $\sigma_i$ s.t. 
$\alpha = \sum_{i} c_{i} \sigma_{i}$, $c_i \in {\bf Z}$.
If $X$ is a simplicial complex, $sk^{i}(X)$ will denote the $i$-skeleton 
of $X$.   
\end{notation}

{\it Proof of Theorem \ref{N2}.}
If we take $A =A_d= \{(1,\epsilon_1,..., \epsilon_d)| 
\epsilon_i \in \{0,1\} 
\}$,        we have that ${\cal I}_{A_d}$ is the ideal of the Segre 
embedding of ${\bf P}^1 \times ... \times {\bf P}^1$ ($d$ times),
i.e. the ideal of the embedding of
 ${\bf P}^1 \times ... \times {\bf P}^1$ ($d$ times)
by the line bundle ${\cal O}(1,..,1)$.
In this case  $\omega=\omega_{d} = (1,0, ..., 0)$ ($0$ repeated $d$ times)
and $n = 2^d$.     
   
Let $b \in {\bf N} A_d$; we have that  $\deg b=(= b \cdot \omega) =k$ if
and only if
$b$ is the sum of $k$ (not necessarily distinct) elements of $A_d$.
By identifying the set $\{1,..., 2^d\}$ with $A_d$, 
we have  that, if $k= \deg b$, 
 $\Delta_b = \cup_{ a_{i_1},...., a_{i_k} \in A_d,
a_{i_1}+...+ a_{i_k} =b} < a_{i_1},...., a_{i_k} >$; we say that 
$< a_{i_1},...., a_{i_k} >$ is a degenerate  $k$-simplex if 
$\exists \; l, m \in \{1, ..., k\}$ with $l \neq m$ 
s.t.  $a_{i_l}= a_{i_m}$; thus $\Delta_b$ 
 is equal to the union of the (possibly degenerate)
$k$-simplexes $S$ with vertices in $A_{d}$
 such that the sum of the vertices (with multiplicities) of $S$ is $b$.
   
 By Theorem \ref{CPS}, in order to prove that  
 ${\cal O}_{{\bf P}^1 \times ... \times {\bf P}^1}(1,..,1)$ 
($d$ times) satisfies $N_2$, we have to prove 
that $H_{1}(\Delta_b)=0$ for each $b \in {\bf N}A_d$ with 
$\deg b  \geq 4$.
Analogously in order to prove that  
 ${\cal O}_{{\bf P}^1 \times ... \times {\bf P}^1}(1,..,1)$ 
($d$ times) satisfies $N_3$, we have  to prove 
that $H_{2}(\Delta_b)=0$ for each $b \in {\bf N}A_d$ 
with $\deg b   \geq 5$.

The proof is by induction on $d$.
Observe that any $b' \in {\bf N} A_{d+1}$ with $\deg b' 
=k $ is equal to 
$\footnotesize{
\left( \hspace{-0.2cm} \begin{array}{c}
b\\
\epsilon
\end{array}
 \hspace{-0.2cm}
\right)}$ for some $b \in {\bf N} A_{d}$ with $\deg b =k$ and 
for some $\varepsilon \in \{0,1,...,k \}$.
Then, in order to prove $N_2$  we suppose (by induction) 
that $H_{1}(\Delta_{b})=0$
$\forall b \in {\bf N}A_d$ with $\deg b  =k$, $k \geq 4$ and 
we  show that
$H_{1}(
\Delta_{
\tiny{
\left( \hspace{-0.2cm} \begin{array}{c}
b\\
\varepsilon
\end{array}
 \hspace{-0.2cm}
\right)}
}
)=0$ for $ \varepsilon \in \{0,...,k\}$ 
 and in order to prove $N_3$
we suppose (by induction) that $H_{2}(\Delta_{b})=0$ 
$ \forall b \in {\bf N} A_d$ with $\deg b =k$,  $k \geq 5 $, and 
we  show that
$H_{2}(
\Delta_{
\tiny{
\left( \hspace{-0.2cm} \begin{array}{c}
b\\
\varepsilon
\end{array}
 \hspace{-0.2cm}
\right)}
})=0$ for  $\varepsilon \in \{0, ..., k\}$.

Observe that, if $\varepsilon \in \{0, k\} $ ($k :=\deg b$), 
 then obviously 
$\Delta_{
\tiny{
\left( \hspace{-0.2cm} \begin{array}{c}
b\\
\varepsilon
\end{array}
 \hspace{-0.2cm}
\right)}
}$ and $\Delta_b$ are isomorphic; besides 
$\Delta_{
\tiny{
\left( \hspace{-0.2cm} \begin{array}{c}
b\\
k- \varepsilon
\end{array}
 \hspace{-0.2cm}
\right)}
}$ is isomorphic to $\Delta_{
\tiny{
\left( \hspace{-0.2cm} \begin{array}{c}
b\\
\varepsilon
\end{array}
 \hspace{-0.2cm}
\right)}
}$ (the isomorphism is given by substituting $0$ with $1$ and $1$ with $0$
in the last coordinate). 
Thus we may consider only the cases $\varepsilon \in \{1,...,[\frac{k}{2}]\}$.
    
First we need some preliminary notation and lemmas.

\begin{notation}
Let $S = \, <a_1, ...., a_k>$ be a (possibly degenerate)
 $k$-simplex, $a_i \in A_d$. Let $\varepsilon \in \{0,...,k\}$.
We denote $$S_{\varepsilon}'=
\cup_{ (\chi_{1},..., \chi_{k}) \; s.t. \; \chi_j \in \{0,1\} \; for
\; j=1,...,k  \; and \; exactly \; \varepsilon \;
of\; \chi_{1},..., \chi_{k} \; are \;equal \;to \;1 }
<
\footnotesize{
\left( \hspace{-0.2cm} \begin{array}{c}
a_{1}\\
\chi_1
\end{array}
 \hspace{-0.2cm}
\right)},
....,
\footnotesize{
\left( \hspace{-0.2cm} \begin{array}{c}
a_{k}\\
\chi_k
\end{array}
 \hspace{-0.2cm}
\right)}
>.$$
\end{notation}

\begin{example}
Let  $S =<a_1, a_2, a_3,a_4>$ is a (possibly degenerate) 
tetrahedron, $a_i \in A_d$. The set 
$S_{1}'$ is the union of the four (possibly degenerate)
 tetrahedrons
$<
\footnotesize{
\left( \hspace{-0.2cm} \begin{array}{c}
a_{1}\\
1
\end{array}
 \hspace{-0.2cm}
\right)},\footnotesize{
\left( \hspace{-0.2cm} \begin{array}{c}
a_{2}\\
0
\end{array}
 \hspace{-0.2cm}
\right)},\footnotesize{
\left( \hspace{-0.2cm} \begin{array}{c}
a_{3}\\
0
\end{array}
 \hspace{-0.2cm}
\right)},
\footnotesize{
\left( \hspace{-0.2cm} \begin{array}{c}
a_{4}\\
0
\end{array}
 \hspace{-0.2cm}
\right)}
>$,
  
$<
\footnotesize{
\left( \hspace{-0.2cm} \begin{array}{c}
a_{1}\\
0
\end{array}
 \hspace{-0.2cm}
\right)},\footnotesize{
\left( \hspace{-0.2cm} \begin{array}{c}
a_{2}\\
1
\end{array}
 \hspace{-0.2cm}
\right)},\footnotesize{
\left( \hspace{-0.2cm} \begin{array}{c}
a_{3}\\
0
\end{array}
 \hspace{-0.2cm}
\right)},
\footnotesize{
\left( \hspace{-0.2cm} \begin{array}{c}
a_{4}\\
0
\end{array}
 \hspace{-0.2cm}
\right)}
>$,
$<
\footnotesize{
\left( \hspace{-0.2cm} \begin{array}{c}
a_{1}\\
0
\end{array}
 \hspace{-0.2cm}
\right)},\footnotesize{
\left( \hspace{-0.2cm} \begin{array}{c}
a_{2}\\
0
\end{array}
 \hspace{-0.2cm}
\right)},\footnotesize{
\left( \hspace{-0.2cm} \begin{array}{c}
a_{3}\\
1
\end{array}
 \hspace{-0.2cm}
\right)},
\footnotesize{
\left( \hspace{-0.2cm} \begin{array}{c}
a_{4}\\
0
\end{array}
 \hspace{-0.2cm}
\right)}
>$,
$<
\footnotesize{
\left( \hspace{-0.2cm} \begin{array}{c}
a_{1}\\
0
\end{array}
 \hspace{-0.2cm}
\right)},\footnotesize{
\left( \hspace{-0.2cm} \begin{array}{c}
a_{2}\\
0
\end{array}
 \hspace{-0.2cm}
\right)},\footnotesize{
\left( \hspace{-0.2cm} \begin{array}{c}
a_{3}\\
0
\end{array}
 \hspace{-0.2cm}
\right)},
\footnotesize{
\left( \hspace{-0.2cm} \begin{array}{c}
a_{4}\\
1
\end{array}
 \hspace{-0.2cm}
\right)}
>$. 
Thus $S_{1}'$ can be obtained from $S$ by ``constructing 
a tetrahedron on everyone
 of the four faces of $S$'' and considering the union
 of these four tetrahedrons.
 The set  $S_{2}'$ is 
 the union of the following six (possibly degenerate) tetrahedrons:
        
$<
\footnotesize{
\left( \hspace{-0.2cm} \begin{array}{c}
a_{1}\\
0
\end{array}
 \hspace{-0.2cm}
\right)},\footnotesize{
\left( \hspace{-0.2cm} \begin{array}{c}
a_{2}\\
0
\end{array}
 \hspace{-0.2cm}
\right)},\footnotesize{
\left( \hspace{-0.2cm} \begin{array}{c}
a_{3}\\
1
\end{array}
 \hspace{-0.2cm}
\right)},
\footnotesize{
\left( \hspace{-0.2cm} \begin{array}{c}
a_{4}\\
1
\end{array}
 \hspace{-0.2cm}
\right)}
>$,
$<
\footnotesize{
\left( \hspace{-0.2cm} \begin{array}{c}
a_{1}\\
0
\end{array}
 \hspace{-0.2cm}
\right)},\footnotesize{
\left( \hspace{-0.2cm} \begin{array}{c}
a_{2}\\
1
\end{array}
 \hspace{-0.2cm}
\right)},\footnotesize{
\left( \hspace{-0.2cm} \begin{array}{c}
a_{3}\\
0
\end{array}
 \hspace{-0.2cm}
\right)},
\footnotesize{
\left( \hspace{-0.2cm} \begin{array}{c}
a_{4}\\
1
\end{array}
 \hspace{-0.2cm}
\right)}
>$,
$<
\footnotesize{
\left( \hspace{-0.2cm} \begin{array}{c}
a_{1}\\
0
\end{array}
 \hspace{-0.2cm}
\right)},\footnotesize{
\left( \hspace{-0.2cm} \begin{array}{c}
a_{2}\\
1
\end{array}
 \hspace{-0.2cm}
\right)},\footnotesize{
\left( \hspace{-0.2cm} \begin{array}{c}
a_{3}\\
1
\end{array}
 \hspace{-0.2cm}
\right)},
\footnotesize{
\left( \hspace{-0.2cm} \begin{array}{c}
a_{4}\\
0
\end{array}
 \hspace{-0.2cm}
\right)}
>$,
          
$<
\footnotesize{
\left( \hspace{-0.2cm} \begin{array}{c}
a_{1}\\
1
\end{array}
 \hspace{-0.2cm}
\right)},\footnotesize{
\left( \hspace{-0.2cm} \begin{array}{c}
a_{2}\\
0
\end{array}
 \hspace{-0.2cm}
\right)},\footnotesize{
\left( \hspace{-0.2cm} \begin{array}{c}
a_{3}\\
0
\end{array}
 \hspace{-0.2cm}
\right)},
\footnotesize{
\left( \hspace{-0.2cm} \begin{array}{c}
a_{4}\\
1
\end{array}
 \hspace{-0.2cm}
\right)}
>$, 
$<
\footnotesize{
\left( \hspace{-0.2cm} \begin{array}{c}
a_{1}\\
1
\end{array}
 \hspace{-0.2cm}
\right)},\footnotesize{
\left( \hspace{-0.2cm} \begin{array}{c}
a_{2}\\
0
\end{array}
 \hspace{-0.2cm}
\right)},\footnotesize{
\left( \hspace{-0.2cm} \begin{array}{c}
a_{3}\\
1
\end{array}
 \hspace{-0.2cm}
\right)},
\footnotesize{
\left( \hspace{-0.2cm} \begin{array}{c}
a_{4}\\
0
\end{array}
 \hspace{-0.2cm}
\right)}
>$, 
$<
\footnotesize{
\left( \hspace{-0.2cm} \begin{array}{c}
a_{1}\\
1
\end{array}
 \hspace{-0.2cm}
\right)},\footnotesize{
\left( \hspace{-0.2cm} \begin{array}{c}
a_{2}\\
1
\end{array}
 \hspace{-0.2cm}
\right)},\footnotesize{
\left( \hspace{-0.2cm} \begin{array}{c}
a_{3}\\
0
\end{array}
 \hspace{-0.2cm}
\right)},
\footnotesize{
\left( \hspace{-0.2cm} \begin{array}{c}
a_{4}\\
0
\end{array}
 \hspace{-0.2cm}
\right)}
>$.
      
Then $S_{2}'$ can be obtained from $S$ by ``constructing 
a tetrahedron on everyone of the six edges of $S$'' and considering 
the union of these six tetrahedrons (see Fig. 1, representing 
$S_{\varepsilon}'$
in the case $S$ is  not degenerate).
\end{example}

\includegraphics[scale=0.53]{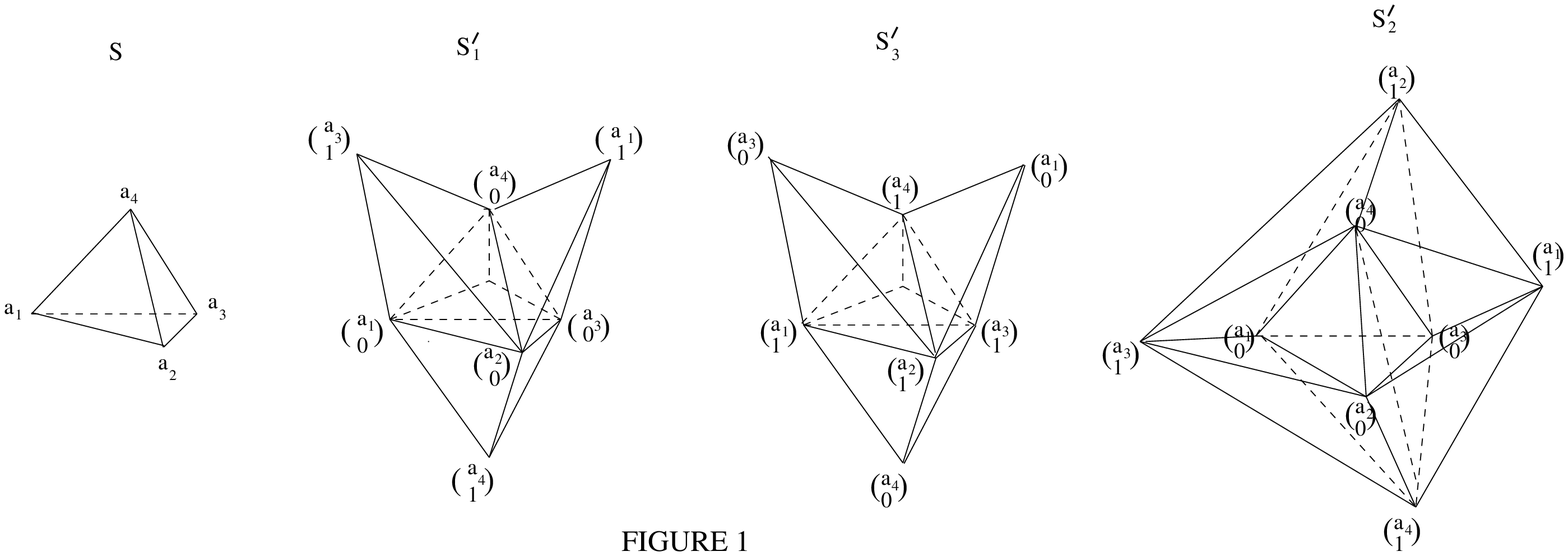}

\vspace{0.15cm}

{\small
{\it NOTE TO FIG. 1.
In the representation of $S_2'$, for the 
sake of simplicity,
we do not represent the tetrahedrons 
$<
\footnotesize{
\left( \hspace{-0.2cm} \begin{array}{c}
a_{1}\\
0
\end{array}
 \hspace{-0.2cm}
\right)},\footnotesize{
\left( \hspace{-0.2cm} \begin{array}{c}
a_{2}\\
1
\end{array}
 \hspace{-0.2cm}
\right)},\footnotesize{
\left( \hspace{-0.2cm} \begin{array}{c}
a_{3}\\
0
\end{array}
 \hspace{-0.2cm}
\right)},
\footnotesize{
\left( \hspace{-0.2cm} \begin{array}{c}
a_{4}\\
1
\end{array}
 \hspace{-0.2cm}
\right)}
>$ and 
$<
\footnotesize{
\left( \hspace{-0.2cm} \begin{array}{c}
a_{1}\\
1
\end{array}
 \hspace{-0.2cm}
\right)},\footnotesize{
\left( \hspace{-0.2cm} \begin{array}{c}
a_{2}\\
0
\end{array}
 \hspace{-0.2cm}
\right)},\footnotesize{
\left( \hspace{-0.2cm} \begin{array}{c}
a_{3}\\
1
\end{array}
 \hspace{-0.2cm}
\right)},
\footnotesize{
\left( \hspace{-0.2cm} \begin{array}{c}
a_{4}\\
0
\end{array}
 \hspace{-0.2cm}
\right)}
>$.}}

\vspace{0.8cm}

Let $b \in {\bf N}A_d$ with $\deg b=k$ and $\varepsilon \in \{0, ..., k\}$.
Obviously
$$\Delta_{
\tiny{
\left( \hspace{-0.2cm} \begin{array}{c}
b\\
\varepsilon
\end{array}
 \hspace{-0.2cm}
\right)}
}= \cup_{S=<a_1,...,a_k>\; with
 \; a_1+...+a_k=b}\; S_{\varepsilon}'.$$

\begin{notation}
Let $b \in {\bf N}A_d$ with $\deg b=k$.
For  $l \in {\bf N}$, $0 \leq l \leq k-1$,
let 
$$F^{l}(\Delta_{b}) =
\cup_{a_1, ..., a_{k} \in A_d \; s.t. \; a_{1}+ .... + a_{k} =b} \; 
\cup_{ i_0,...,i_l \in \{1,...., k \} }  
<
\footnotesize{
\left( \hspace{-0.2cm} \begin{array}{c}
a_{i_0}\\
0
\end{array}
 \hspace{-0.2cm}
\right)},
...,
\footnotesize{
\left( \hspace{-0.2cm} \begin{array}{c}
a_{i_l}\\
0
\end{array}
 \hspace{-0.2cm}
\right)}
>.$$ 
\end{notation}


Observe that  $F^{l}(\Delta_b)  \subseteq \Delta_{
\tiny{
\left( \hspace{-0.2cm} \begin{array}{c}
b\\
\varepsilon
\end{array}
 \hspace{-0.2cm}
\right)}
}$ iff $k- \varepsilon \geq l+1$.

The idea of the proof is to consider a $l$-cycle (for $l=1,2$) in 
$\Delta_{
\tiny{
\left( \hspace{-0.2cm} \begin{array}{c}
b\\
\varepsilon
\end{array}
 \hspace{-0.2cm}
\right)}
}$  and to show that it is homologous to a $l$-cycle
in $F^{l}(\Delta_b)$  and then
 to show that it is homologous to $0$ by 
using that $H_{l}(\Delta_b)=0$.

\begin{remark} \label{cone}
Let  $S= <a_1,..., a_k>$, $a_i \in A_d$.
If $k \geq 4$ and $1 \leq  \varepsilon \leq k-2$,  
the set $S_{\varepsilon}'$ contains 
the  cone with vertex 
$\footnotesize{
\left( \hspace{-0.2cm} \begin{array}{c}
a_{l}\\
1
\end{array}
 \hspace{-0.2cm}
\right)}$  on the border of  
$<
\footnotesize{
\left( \hspace{-0.2cm} \begin{array}{c}
a_{i_1}\\
0
\end{array}
 \hspace{-0.2cm}
\right)},
\footnotesize{
\left( \hspace{-0.2cm} \begin{array}{c}
a_{i_{2}}\\
0
\end{array}
 \hspace{-0.2cm}
\right)},
\footnotesize{
\left( \hspace{-0.2cm} \begin{array}{c}
a_{i_{3}}\\
0
\end{array}
 \hspace{-0.2cm}
\right)}
>$ for any $i_1, i_2, i_{3} ,l \in \{1,...,k\}$ with $l \neq i_j$
for $j=1,2,3$. This is true in particular if $k \geq 4$ and $\varepsilon 
\in \{1,...,[\frac{k}{2}]\}$.

If $k \geq 5 $ and $1 \leq \varepsilon \leq k-3$,
the set $S_{\varepsilon }'$ contains 
the  cone with vertex 
$\footnotesize{
\left( \hspace{-0.2cm} \begin{array}{c}
a_{l}\\
1
\end{array}
 \hspace{-0.2cm}
\right)}$  on the border of 
$< \footnotesize{
\left( \hspace{-0.2cm} \begin{array}{c}
a_{i_1}\\
0
\end{array}
 \hspace{-0.2cm}
\right)},...,
\footnotesize{
\left( \hspace{-0.2cm} \begin{array}{c}
a_{i_4}\\
0
\end{array}
 \hspace{-0.2cm}
\right)}>$ for any $i_1,...., i_4 ,l \in \{1,..., k\}$ with $l \neq i_j$
for $j=1,2,3,4$.  
This is true in particular if $k \geq 5$ and 
$\varepsilon \in \{1,...,[\frac{k}{2}]\}$.
\end{remark}

\begin{definition} \label{R}
For any $c \in {\bf N} A_d$ with $\deg c =s$ and $\varepsilon \in 
\{1,..., s\} $, we define $R_{c, \varepsilon}$ the following  set: 
$$ \cup_{\alpha_{1}, ..., \alpha_{s} \in A_{d} \; s.t. \;
 \alpha_{1}+ ...+ \alpha_{s}=c }  \;
\cup_{i_{1},..., i_{s-1} \in \{1,..., s \}, i_{l}\neq i_{m}}
<\footnotesize{
\left( \hspace{-0.2cm} 
\begin{array}{c}
\alpha_{i_{1}}\\
1
\end{array}
 \hspace{-0.2cm}
\right)
}, ...,
\footnotesize{
\left( \hspace{-0.2cm} 
\begin{array}{c}
\alpha_{i_{\varepsilon -1}}\\
1
\end{array}
 \hspace{-0.2cm}
\right)
}, 
\footnotesize{
\left( \hspace{-0.2cm} 
\begin{array}{c}
\alpha_{i_{\varepsilon}}\\
0
\end{array}
 \hspace{-0.2cm}
\right)
},...,   \footnotesize{
\left( \hspace{-0.2cm} 
\begin{array}{c}
\alpha_{i_{s-1}}\\
0
\end{array}
 \hspace{-0.2cm}
\right)
}     >. $$
\end{definition}

\begin{lemma} \label{Rcon}
Let $c \in {\bf N} A_d$ with $\deg c =s$.
We have that $\tilde{H}_{i}(
\Delta_{
\tiny{
\left( \hspace{-0.2cm} \begin{array}{c}
c\\
\varepsilon -1
\end{array}
 \hspace{-0.2cm}
\right)}
}) = 0$ implies
$\tilde{H}_{i}(R_{c, \varepsilon})=0$ 
if we are in one of the following cases:
a) $i=0$, $s \geq 3$, $\varepsilon \in \{1,...,[\frac{s+1}{2} ]\}$  
b) $i=1$, $s \geq 4$, $\varepsilon \in \{1,...,[\frac{s+1}{2} ]\}$  
\end{lemma}

{\it Proof.}
Observe that $R_{c, \varepsilon} \subseteq \Delta_{
\tiny{
\left( \hspace{-0.2cm} \begin{array}{c}
c\\
\varepsilon -1
\end{array}
 \hspace{-0.2cm}
\right)}}$.
Since $\tilde{H}_{i}(\Delta_{
\tiny{
\left( \hspace{-0.2cm} \begin{array}{c}
c\\
\varepsilon -1
\end{array}
 \hspace{-0.2cm}
\right)}
}) = 0$, we have  $\tilde{H}_{i}(sk^{i+1}
(\Delta_{
\tiny{
\left( \hspace{-0.2cm} \begin{array}{c}
c\\
\varepsilon -1
\end{array}
 \hspace{-0.2cm}
\right)}
})) = 0$. 
Obviously $sk^{i+1}(R_{c, \varepsilon}) \subseteq sk^{i+1}
(\Delta_{
\tiny{
\left( \hspace{-0.2cm} \begin{array}{c}
c\\
\varepsilon -1
\end{array}
 \hspace{-0.2cm}
\right)}
})$. We want to show $\tilde{H}_{i}(sk^{i+1}(R_{c, \varepsilon}))=0$. 
Let $\beta $ be a $i$-cycle in $sk^{i+1}(R_{c, \varepsilon})$. Since  
 $\tilde{H}_{i}(sk^{i+1}
(\Delta_{
\tiny{
\left( \hspace{-0.2cm} \begin{array}{c}
c\\
\varepsilon -1
\end{array}
 \hspace{-0.2cm}
\right)}
})) = 0$, there exists  a $(i+1)$-chain $\eta$ in $sk^{i+1}
(\Delta_{
\tiny{
\left( \hspace{-0.2cm} \begin{array}{c}
c\\
\varepsilon -1
\end{array}
 \hspace{-0.2cm}
\right)}
})$ s.t. $\partial \eta = \beta$. Suppose $sp(\eta ) = \cup_{j} F_j $,
 where
$F_j$ are $(i+1)$-simplexes in    $sk^{i+1}
(\Delta_{
\tiny{
\left( \hspace{-0.2cm} \begin{array}{c}
c\\
\varepsilon -1
\end{array}
 \hspace{-0.2cm}
\right)}
})$; consider now a $(i+1)$-chain $\psi$  in  
$sk^{i+1}(R_{c, \varepsilon})$
whose support is $\cup_j \hat{F}_j$, where $\hat{F}_j = F_j$ if $F_j
\subseteq  sk^{i+1}(R_{c, \varepsilon})$ and $\hat{F}_{j}$
is a cone on the border of $F_j$
if $F_j \not\subseteq sk^{i+1}(R_{c, \varepsilon})$, in such way that 
$\partial \psi = \beta $ (observe that in our cases 
such cones exist, in fact:
$R_{c, \varepsilon } $ is the union of the (possibly degenerate)
$(s-1)$-simplexes
``obtained from the (possibly degenerate) $s$-simplexes of $\Delta_{
\tiny{
\left( \hspace{-0.2cm} \begin{array}{c}
c\\
\varepsilon -1
\end{array}
 \hspace{-0.2cm}
\right)}}$ by taking off  a vertex whose last coordinate is $0$'';
 in the case 
$i=0$ one can check that the $1$-simplexes whose vertices have the last 
coordinates equal to $1$, $1$ or to $1$,  $0$ are contained in $R_{c,
 \varepsilon}$, while for a  $1$-simplex $F$ whose vertices have the last 
coordinates equal to $0$, $0$, there exists a cone, $\hat{F}$, on the 
 border of $F$
with $\hat{F} \subseteq R_{c, \varepsilon}$, since $s \geq 3$; analogously
the case b)). Thus we proved  $\tilde{H}_{i}(sk^{i+1}(R_{c, 
\varepsilon}))=0$.
Thus  $\tilde{H}_{i}(R_{c, \varepsilon})=0$.
        
\hfill {\it Q.e.d. in Lemma \ref{Rcon} }

\vspace{0.2cm}

\begin{center}
PROOF THAT ${\cal O}_{{\bf P}^{1} 
\times ...\times {\bf P}^{1} }(1,..., 1)$ SATISFIES PROPERTY $N_2$.
\end{center}

\begin{lemma} \label{sc}
Let $b \in {\bf N} A_d$, $\deg b =k$,  $k
\geq  4$ and $ \varepsilon \in \{1,...,
[\frac{k}{2}]\}$.  
Every $1-$cycle $\gamma$ in 
$\Delta_{
\tiny{
\left( \hspace{-0.2cm} \begin{array}{c}
b\\
\varepsilon
\end{array}
 \hspace{-0.2cm}
\right)}
}$ is homologous to a $1-$cycle in 
$F^{1}(\Delta_{b})$ (which is $\subseteq \Delta_{
\tiny{
\left( \hspace{-0.2cm} \begin{array}{c}
b\\
\varepsilon
\end{array}
 \hspace{-0.2cm}
\right)}
}$ since $k - \varepsilon \geq 2$). 
\end{lemma}

{\it Proof.}
Obviously we can suppose $sp(\gamma) \subseteq sk^{1}(\Delta_{
\tiny{\left( \hspace{-0.2cm} \begin{array}{c}
b\\
\varepsilon
\end{array}
 \hspace{-0.2cm}
\right)}
})$. 
 The proof is by induction on the cardinality of  
$(sp(\gamma) \cap sk^{0}(\Delta_{
\tiny{
\left( \hspace{-0.2cm} \begin{array}{c}
b\\
\varepsilon
\end{array}
 \hspace{-0.2cm}
\right)}
})) -F^{1}(\Delta_b)$, i.e. we will prove that $\gamma$ is homotopically 
equivalent a $1$-cycle
 $\tilde{\gamma}$ s.t. $\sharp( (sp(\tilde{\gamma}) \cap  sk^{0}(\Delta_{
\tiny{
\left( \hspace{-0.2cm} \begin{array}{c}
b\\
\varepsilon
\end{array}
 \hspace{-0.2cm}
\right)}
})) - F^{1}(\Delta_b)) <
\sharp ((sp(\gamma) \cap  sk^{0}(\Delta_{
\tiny{
\left( \hspace{-0.2cm} \begin{array}{c}
b\\
\varepsilon
\end{array}
 \hspace{-0.2cm}
\right)}
})) - F^{1}(\Delta_b))$. 
        
Let $ \footnotesize{
\left( \hspace{-0.2cm} \begin{array}{c}
a\\
1
\end{array}
 \hspace{-0.2cm}
\right)} \in 
(sp(\gamma) \cap sk^{0}(\Delta_{
\tiny{
\left( \hspace{-0.2cm} \begin{array}{c}
b\\
\varepsilon
\end{array}
 \hspace{-0.2cm}
\right)}
}))  
- F^{1}(\Delta_b) $,  ($a \in A_d$). 
Let $< \footnotesize{ \left( \hspace{-0.2cm} \begin{array}{c}
a\\
1
\end{array}
 \hspace{-0.2cm}
\right)}, P_{1}> \cup 
<\footnotesize{
\left( \hspace{-0.2cm} \begin{array}{c}
a\\
1
\end{array}
 \hspace{-0.2cm}
\right)} , P_{2}>  \; \subseteq sp(\gamma)$, 
with  $P_j \in sk^{0}( 
\Delta_{
\tiny{
\left( \hspace{-0.2cm} \begin{array}{c}
b\\
\varepsilon
\end{array}
 \hspace{-0.2cm}
\right)}})$ for $j=1,2$.
 Precisely let $\gamma = \sigma_1 + \sigma_2 + ...$, where
$\sigma_1$ and $\sigma_2 $ are two simplexes 
 $\sigma_1: [0,1] \rightarrow \;
 < \footnotesize{ \left( \hspace{-0.2cm} \begin{array}{c}
a\\
1
\end{array}
 \hspace{-0.2cm}
\right)}, P_{1}>$ and $ \sigma_2: [0,1] \rightarrow \;
<\footnotesize{
\left( \hspace{-0.2cm} \begin{array}{c}
a\\
1
\end{array}
 \hspace{-0.2cm}
\right)} , P_{2}> $ s.t. $\sigma_{1} (0)= P_1$, $\sigma_{1} (1)
 = \footnotesize{ \left(
 \hspace{-0.2cm} \begin{array}{c}
a\\
1
\end{array}
 \hspace{-0.2cm}
\right)} = \sigma_2 (0)$, $\sigma_{2}(1) = P_2$.

Let $\alpha$ be the $1$-cycle  $ - \sigma_1  - \sigma_2 + \sigma'_{1}
+ \sigma'_{2}$, where 
$\sigma'_1$ and $\sigma'_2 $ are two simplexes 
 $\sigma'_1: [0,1] \rightarrow \;
 < \footnotesize{ \left( \hspace{-0.2cm} \begin{array}{c}
a\\
0
\end{array}
 \hspace{-0.2cm}
\right)}, P_{1}>$ and $ \sigma'_2: [0,1] \rightarrow \;
<\footnotesize{
\left( \hspace{-0.2cm} \begin{array}{c}
a\\
0
\end{array}
 \hspace{-0.2cm}
\right)} , P_{2}> $ s.t. $\sigma'_{1} (0)= P_1$, $\sigma'_{1} (1)
 = \footnotesize{ \left(
 \hspace{-0.2cm} \begin{array}{c}
a\\
0
\end{array}
 \hspace{-0.2cm}
\right)}= \sigma'_2 (0)$, $\sigma'_{2} (1)= P_{2}$.

The  support of $\alpha$  is the union of the two cones 
with vertices $\footnotesize{
\left( \hspace{-0.2cm} \begin{array}{c}
a\\
1
\end{array}
 \hspace{-0.2cm}
\right)}$ and 
$ \footnotesize{
\left( \hspace{-0.2cm} \begin{array}{c}
a\\
0
\end{array}
 \hspace{-0.2cm}
\right)}$ on  $\{P_1, P_2\}$. 
 
We state that  $P_i \in R_{b-a, \varepsilon}$ for $i = 1,2$. 
In fact, $<P_i, \footnotesize{
\left( \hspace{-0.2cm} \begin{array}{c}
a\\
1
\end{array}
 \hspace{-0.2cm}
\right)}> \;  \subseteq \Delta_{
\tiny{
\left( \hspace{-0.2cm} \begin{array}{c}
b\\
\varepsilon
\end{array}
 \hspace{-0.2cm}
\right)}} $, then 
$ P_i  
\in
\Delta_{
\tiny{
\left( \hspace{-0.2cm} \begin{array}{c}
b -a\\
\varepsilon -1
\end{array}
 \hspace{-0.2cm}
\right)}} $; we recall that $R_{b-a, \varepsilon }$ is 
$$\cup_{\alpha_{1}, ..., \alpha_{k-1} \in A_{d} \; s.t. \;
 \alpha_{1}+ ...+ \alpha_{k-1}=b-a }  \;
\cup_{i_{1},..., i_{k-2} \in \{1,..., k-1 \}, i_{l}\neq i_{m}}
<\footnotesize{
\left( \hspace{-0.2cm} 
\begin{array}{c}
\alpha_{i_{1}}\\
1
\end{array}
 \hspace{-0.2cm}
\right)
}, ...,
\footnotesize{
\left( \hspace{-0.2cm} 
\begin{array}{c}
\alpha_{i_{\varepsilon -1}}\\
1
\end{array}
 \hspace{-0.2cm}
\right)
}, 
\footnotesize{
\left( \hspace{-0.2cm} 
\begin{array}{c}
\alpha_{i_{\varepsilon}}\\
0
\end{array}
 \hspace{-0.2cm}
\right)
},...,   \footnotesize{
\left( \hspace{-0.2cm} 
\begin{array}{c}
\alpha_{i_{k-2}}\\
0
\end{array}
 \hspace{-0.2cm}
\right)
}     >,$$
i.e. $R_{b-a, \varepsilon }$ is the 
union of the (possibly degenerate) $(k-2)$-simplexes
``obtained from the (possibly degenerate) $(k-1)$-simplexes of 
 $\Delta_{
\tiny{
\left( \hspace{-0.2cm} \begin{array}{c}
b -a \\
\varepsilon -1
\end{array}
 \hspace{-0.2cm}
\right)}} $ by taking off a vertex whose last coordinate is $0$''; then,
if the last coordinate of $P_i$ is $1$, we may conclude at once that 
 $P_{i} \in  R_{b-a, \varepsilon}$; also if 
 the last coordinate of $P_i$ is $0$, we may conclude  that 
 $P_{i} \in  R_{b-a, \varepsilon}$,
because 
the number of the vertices whose last coordinate is $0$ in a
(possibly degenerate) $(k-1)$-simplex of 
 $\Delta_{
\tiny{
\left( \hspace{-0.2cm} \begin{array}{c}
b -a \\
\varepsilon -1
\end{array}
 \hspace{-0.2cm}
\right)}} $ is
$k-1 - (\varepsilon -1) \geq 2$.

Thus $sp(\alpha) \subseteq  C$, where $C$ is the union of the two
 cones $ < \footnotesize{
\left( \hspace{-0.2cm} \begin{array}{c}
a\\
1
\end{array}
 \hspace{-0.2cm}
\right)}, R_{b-a, \varepsilon }>$ and   
$
< \footnotesize{
\left( \hspace{-0.2cm} \begin{array}{c}
a\\
0
\end{array}
 \hspace{-0.2cm}
\right)}, R_{b-a, \varepsilon }>$.
Observe that  $ C \subseteq
\Delta_{
\tiny{
\left( \hspace{-0.2cm} \begin{array}{c}
b\\
\varepsilon
\end{array}
 \hspace{-0.2cm}
\right)}}$.
  
We want to show that $H_{1}(C)=0$:
by Theorem \ref{CPS}, since  ${\cal O}_{{\bf P}^{1} \times ...\times 
{\bf P}^{1}}(1,..., 1)$ ($d$ times) satisfies
Property $N_1$ $\forall d$, we have  $\tilde{H}_{0}(\Delta_g)=0 $ 
$ \forall g \in {\bf N}A_d$ with $\deg g \geq 3 $, $\forall d$
(this can be easily proved directly without using   
that  ${\cal O}_{{\bf P}^{1} \times ...\times 
{\bf P}^{1}}(1,..., 1)$  satisfies
Property $N_1$);
then $\tilde{H}_{0}(
\Delta_{
\tiny{
\left( \hspace{-0.2cm} \begin{array}{c}
b-a\\
\varepsilon -1
\end{array}
 \hspace{-0.2cm}
\right)}
}) = 0$; thus $\tilde{ H}_{0}(R_{b-a, \varepsilon})=0$ by Lemma \ref{Rcon};
we have $\tilde{H}_{i}(C)=\tilde{H}_{i-1}(R_{b-a, \varepsilon })$; thus 
$H_{1}(C)=\tilde{ H}_{0}(R_{b-a, \varepsilon})=0$.
        
Thus we have that $\alpha$ is homologous to $0$.

Thus   $\gamma $
is homologous to $\gamma +\alpha $; obviously 
$\gamma + \alpha$ (and then $\gamma$)  is homologous 
to a $1$-cycle $\tilde{\gamma}$ whose support can be obtained from 
$sp(\gamma)$ by substituting  the cone with vertex $ \footnotesize{
\left( \hspace{-0.2cm} \begin{array}{c}
a\\
1
\end{array}
 \hspace{-0.2cm}
\right)}$ on  $\{P_1, P_2\}$
with the cone with vertex $\footnotesize{
\left( \hspace{-0.2cm} \begin{array}{c}
a\\
0
\end{array}
 \hspace{-0.2cm}
\right)}$ on  $\{P_1, P_2 \}$.
 Then $\sharp( (sp(\tilde{\gamma}) \cap  sk^{0}(\Delta_{
\tiny{
\left( \hspace{-0.2cm} \begin{array}{c}
b\\
\varepsilon
\end{array}
 \hspace{-0.2cm}
\right)}
})) - F^{1}(\Delta_b)) <
\sharp ((sp(\gamma) \cap  sk^{0}(\Delta_{
\tiny{
\left( \hspace{-0.2cm} \begin{array}{c}
b\\
\varepsilon
\end{array}
 \hspace{-0.2cm}
\right)}
})) - F^{1}(\Delta_b))$; thus we  conclude.
  \hfill {\it Q.e.d. in Lemma \ref{sc}}

\vspace{0.3cm}

 In order to prove that 
 ${\cal O}_{{\bf P}^1 \times ... \times {\bf P}^1}(1,..,1)$ 
($d$ times) satisfies $N_2$ for any $ d$,
we suppose (by induction) 
that $H_{1}(\Delta_{b})=0$
$\forall b \in {\bf N}A_d$ with $\deg b  =k $, $k \geq 4$ and 
we  show that
$H_{1}(
\Delta_{
\tiny{
\left( \hspace{-0.2cm} \begin{array}{c}
b\\
\varepsilon
\end{array}
 \hspace{-0.2cm}
\right)}
}
)=0$ for $ \varepsilon \in \{1,..., [\frac{k}{2}]\}$.

\underline{Cases $\varepsilon  \leq k-3 $.}
We know that every $1$-cycle 
$\gamma$  in $\Delta_{
\tiny{
\left( \hspace{-0.2cm} \begin{array}{c}
b\\
\varepsilon
\end{array}
 \hspace{-0.2cm}
\right)}
} $ is homologous to a $1-$cycle in 
$F^{1}(\Delta_{b})$ by Lemma \ref{sc}.  Thus, since $F^{2}(\Delta_b) 
\subseteq \Delta_{
\tiny{
\left( \hspace{-0.2cm} \begin{array}{c}
b\\
\varepsilon
\end{array}
 \hspace{-0.2cm}
\right)}}$ and
$H_{1}(F^{2}(\Delta_b))=0$
(because, by induction hypothesis,
 $H_{1}(\Delta_b)=0$), we have that $H_{1}(
\Delta_{
\tiny{
\left( \hspace{-0.2cm} \begin{array}{c}
b\\
\varepsilon
\end{array}
 \hspace{-0.2cm}
\right)}
})=0$.

\underline{Cases $\varepsilon  > k-3$.}
These cases are slightly more difficult.
By  Lemma \ref{sc} every $1$-cycle $\gamma$ in $\Delta_{\tiny{
\left( \hspace{-0.2cm} \begin{array}{c}
b\\
\varepsilon 
\end{array}
 \hspace{-0.2cm}
\right)}
}$ is homologous to a $1$-cycle
$\gamma'$ in $F^{1}(\Delta_b)$.
But in these cases we have not the inclusion 
$F^{2}(\Delta_b) \subseteq 
\Delta_{\tiny{
\left( \hspace{-0.2cm} \begin{array}{c}
b\\
\varepsilon
\end{array}
 \hspace{-0.2cm}
\right)}
}$, thus we have to conclude the proof in another way.

Since $H_{1} (F^{2}(\Delta_{b})) =0 $,
  there exists a $2$-chain $\mu$  in $F^{2}(\Delta_b)$ s.t.
 $\partial \mu= \gamma'$. Let
$sp(\mu) = \cup_{i} F_i$, $F_i$ triangles in $F^{2}(\Delta_b)$.
Consider a $2$-chain 
$\psi$ in  $\Delta_{\tiny{
\left( \hspace{-0.2cm} \begin{array}{c}
b\\
\varepsilon
\end{array}
 \hspace{-0.2cm}
\right)}
}$ whose support is   
$\cup_{i} \hat{F_i}$, where
$\hat{F_i}$ is a cone $\subseteq \Delta_{\tiny{
\left( \hspace{-0.2cm} \begin{array}{c}
b\\
\varepsilon
\end{array}
 \hspace{-0.2cm}
\right)}
}$ on the border of $F_i$ (there exists by Remark \ref{cone}),
in such way  that $\partial \psi = \gamma'$;
thus $[\gamma'] =0 $ in $H_{1}(\Delta_{\tiny{
\left( \hspace{-0.2cm} \begin{array}{c}
b\\
\varepsilon
\end{array}
 \hspace{-0.2cm}
\right)}
})$, thus $[\gamma] =0 $ in $H_{1}(\Delta_{\tiny{
\left( \hspace{-0.2cm} \begin{array}{c}
b\\
\varepsilon
\end{array}
 \hspace{-0.2cm}
\right)}
})$. Thus $H_{1}(\Delta_{\tiny{
\left( \hspace{-0.2cm} \begin{array}{c}
b\\
\varepsilon
\end{array}
 \hspace{-0.2cm}
\right)}
})=0$.

\vspace{0.2cm}

\begin{center}
PROOF THAT ${\cal O}_{{\bf P}^{1} 
\times ...\times {\bf P}^{1} }(1,..., 1)$ SATISFIES PROPERTY  $N_3$
\end{center}

\begin{lemma} \label{sc1}
Let $b \in {\bf N} A_d$ with $\deg b=k $ and $\varepsilon \in 
\{1,...,[\frac{k}{2}]\}$. 
If $k \geq 5$, every $2$-cycle $\mu $ in 
 $\Delta_{
\tiny{
\left( \hspace{-0.2cm} \begin{array}{c}
b\\
\varepsilon
\end{array}
 \hspace{-0.2cm}
\right)} }$
is homologous    
to a $2$-cycle  in $F^{2}(\Delta_b) $ (which is $\subseteq \Delta_{
\tiny{
\left( \hspace{-0.2cm} \begin{array}{c}
b\\
\varepsilon
\end{array}
 \hspace{-0.2cm}
\right)} }$ since $k - \varepsilon \geq 3$).
\end{lemma}

{\it Proof.}
Obviously we can suppose that  $sp(\mu) \subseteq sk^{2}(\Delta_{
\tiny{\left( \hspace{-0.2cm} \begin{array}{c}
b\\
\varepsilon
\end{array}
 \hspace{-0.2cm}
\right)}
})$. 
The proof is by induction on the cardinality of  
$(sp(\mu) \cap sk^{0}(\Delta_{
\tiny{
\left( \hspace{-0.2cm} \begin{array}{c}
b\\
\varepsilon
\end{array}
 \hspace{-0.2cm}
\right)}
})) -F^{2}(\Delta_b)$, i.e. we will  prove that $\mu $ is homotopically 
equivalent to a $2$-cycle $\tilde{\mu}$ s.t. 
$\sharp ((sp(\tilde{\mu}) \cap sk^{0}(\Delta_{
\tiny{
\left( \hspace{-0.2cm} \begin{array}{c}
b\\
\varepsilon
\end{array}
 \hspace{-0.2cm}
\right)}
})) -F^{2}(\Delta_b) ) < 
\sharp ((sp(\mu) \cap sk^{0}(\Delta_{
\tiny{
\left( \hspace{-0.2cm} \begin{array}{c}
b\\
\varepsilon
\end{array}
 \hspace{-0.2cm}
\right)}
})) -F^{2}(\Delta_b))$. 
        
Let $ \footnotesize{
\left( \hspace{-0.2cm} \begin{array}{c}
a\\
1
\end{array}
 \hspace{-0.2cm}
\right)} \in 
(sp(\mu) \cap sk^{0}(\Delta_{
\tiny{
\left( \hspace{-0.2cm} \begin{array}{c}
b\\
\varepsilon
\end{array}
 \hspace{-0.2cm}
\right)}
}))  
- F^{2}(\Delta_b) $,  ($a \in A_d$). 
Let $< \footnotesize{ \left( \hspace{-0.2cm} \begin{array}{c}
a\\
1
\end{array}
 \hspace{-0.2cm}
\right)}, P_{1}, P_{2}> \cup 
<\footnotesize{
\left( \hspace{-0.2cm} \begin{array}{c}
a\\
1
\end{array}
 \hspace{-0.2cm}
\right)} , P_{2}, P_{3}> \cup ....... \cup 
<\footnotesize{
\left( \hspace{-0.2cm} \begin{array}{c}
a\\
1
\end{array}
 \hspace{-0.2cm}
\right)} , P_{r-1}, P_{r}>
\cup 
<\footnotesize{
\left( \hspace{-0.2cm} \begin{array}{c}
a\\
1
\end{array}
 \hspace{-0.2cm}
\right)} , P_{r}, P_{1}> \; \subseteq sp(\mu)$, 
with  $P_j \in sk^{0}( 
\Delta_{
\tiny{
\left( \hspace{-0.2cm} \begin{array}{c}
b\\
\varepsilon
\end{array}
 \hspace{-0.2cm}
\right)}})$. 
 
Let $\alpha $ be a $2$-cycle whose support is the union of the two cones 
with vertices $\footnotesize{
\left( \hspace{-0.2cm} \begin{array}{c}
a\\
1
\end{array}
 \hspace{-0.2cm}
\right)}$ and 
$ \footnotesize{
\left( \hspace{-0.2cm} \begin{array}{c}
a\\
0
\end{array}
 \hspace{-0.2cm}
\right)}$ on the polygon with vertices $P_1, ..., P_r$; choose $\alpha$
in  such way that
$\mu + \alpha $ is  homologous 
to a $2$-cycle $\tilde{\mu}$ whose support can be obtained from 
$sp(\mu)$ by substituting  the cone with vertex $ \footnotesize{
\left( \hspace{-0.2cm} \begin{array}{c}
a\\
1
\end{array}
 \hspace{-0.2cm}
\right)}$ on the polygon with vertices $P_1, ..., P_r$
with the cone with vertex $\footnotesize{
\left( \hspace{-0.2cm} \begin{array}{c}
a\\
0
\end{array}
 \hspace{-0.2cm}
\right)}$ on the polygon with vertices $P_1, ..., P_r$. 
 
We state that 
$<P_i, P_{i+1}> \; \subseteq R_{b-a, \varepsilon}$ for $i=1,..., r-1$ and 
$<P_r, P_{1}> \; \subseteq R_{b-a, \varepsilon}$.
In fact: 
   
$<P_i, P_{i+1}, \; \footnotesize{
\left( \hspace{-0.2cm} \begin{array}{c}
a\\
1
\end{array}
 \hspace{-0.2cm}
\right)}>  \; \subseteq  \Delta_{
\tiny{
\left( \hspace{-0.2cm} \begin{array}{c}
b\\
\varepsilon
\end{array}
 \hspace{-0.2cm}
\right)} }$, then 
$<P_i, P_{i+1}> 
\; \subseteq  \Delta_{
\tiny{
\left( \hspace{-0.2cm} \begin{array}{c}
b -a \\
\varepsilon -1
\end{array}
 \hspace{-0.2cm}
\right)} }$;
since 
$ R_{b-a, \varepsilon}$ 
 is the union of the (possibly degenerate) $(k-2)$-simplexes
``obtained from the (possibly degenerate) $(k-1)$-simplexes of 
 $\Delta_{
\tiny{
\left( \hspace{-0.2cm} \begin{array}{c}
b -a \\
\varepsilon -1
\end{array}
 \hspace{-0.2cm}
\right)}} $ by taking off a vertex whose last coordinate is $0$'' and since 
the number of the vertices whose last coordinate is $0$ in a
(possibly degenerate) $(k-1)$-simplex of 
 $\Delta_{
\tiny{
\left( \hspace{-0.2cm} \begin{array}{c}
b -a \\
\varepsilon -1
\end{array}
 \hspace{-0.2cm}
\right)}} $ is
$k-1 - (\varepsilon -1) \geq 3$, we have $<P_{i}, P_{i+1}> \;  
\subseteq  R_{b-a, \varepsilon}$ (obviously 
in a completely analogous way we have $<P_r, P_{1}> \; \subseteq
 R_{b-a, \varepsilon}$).

Thus $sp(\alpha) \subseteq  C$, where $C$ is the union of the two
 cones $ < \footnotesize{
\left( \hspace{-0.2cm} \begin{array}{c}
a\\
1
\end{array}
 \hspace{-0.2cm}
\right)}, R_{b-a, \varepsilon }>$ and   
$
< \footnotesize{
\left( \hspace{-0.2cm} \begin{array}{c}
a\\
0
\end{array}
 \hspace{-0.2cm}
\right)}, R_{b-a, \varepsilon }>$.
Observe that  $ C \subseteq
\Delta_{
\tiny{
\left( \hspace{-0.2cm} \begin{array}{c}
b\\
\varepsilon
\end{array}
 \hspace{-0.2cm}
\right)}}$.

We want to show that $H_{2}(C)=0$: 
we have already proved that  ${\cal O}_{{\bf P}^1 
\times ... \times {\bf P}^1}(1,..,1)$  satisfies Property $N_2$ 
i.e. $H_{1}(\Delta_{g})=0 $ $ \forall g$ with $\deg g \geq 4$;
thus
$H_{1}(
\Delta_{
\tiny{
\left( \hspace{-0.2cm} \begin{array}{c}
b-a\\
\varepsilon -1
\end{array}
 \hspace{-0.2cm}
\right)}
}) = 0$; then  $H_{1}(R_{b-a,\varepsilon })=0$ by
 Lemma \ref{Rcon};
we have
$\tilde{H}_{i}(C)=\tilde{H}_{i-1}(R_{b-a, \varepsilon })$; thus 
 $H_{2}(C)= H_{1}(R_{b-a, \varepsilon})=0$.

Thus we have that $\alpha$ is homologous to $0$.

Thus  $\mu  $
is homologous to $\mu + \alpha $; the cycle
$\mu + \alpha$ (and then $\mu$) is homologous 
to a $2$-cycle $\tilde{\mu}$ whose support can be obtained from 
$sp(\mu)$ by substituting  the cone with vertex $ \footnotesize{
\left( \hspace{-0.2cm} \begin{array}{c}
a\\
1
\end{array}
 \hspace{-0.2cm}
\right)}$ on the polygon with vertices $P_1, ..., P_r$
with the cone with vertex $\footnotesize{
\left( \hspace{-0.2cm} \begin{array}{c}
a\\
0
\end{array}
 \hspace{-0.2cm}
\right)}$ on the polygon with vertices $P_1, ..., P_r$.
 Then  $\sharp ((sp(\tilde{\mu}) \cap  sk^{0}(\Delta_{
\tiny{
\left( \hspace{-0.2cm} \begin{array}{c}
b\\
\varepsilon
\end{array}
 \hspace{-0.2cm}
\right)}
})) - F^{2}(\Delta_b)) <
\sharp ((sp(\mu) \cap  sk^{0}(\Delta_{
\tiny{
\left( \hspace{-0.2cm} \begin{array}{c}
b\\
\varepsilon
\end{array}
 \hspace{-0.2cm}
\right)}
})) - F^{2}(\Delta_b))$; thus we  conclude.
  \hfill{\it Q.e.d. in Lemma \ref{sc1}}

\vspace{0.3cm}

 In order to prove that 
 ${\cal O}_{{\bf P}^1 \times ... \times {\bf P}^1}(1,..,1)$ 
($d$ times) satisfies 
 $N_3$ for any $d$,  we suppose (by induction) 
that $H_{2}(\Delta_{b})=0$
$\forall b \in {\bf N}A_d$ with $\deg b  =k$, $k \geq 5$ and 
we  show that
$H_{2}(
\Delta_{
\tiny{
\left( \hspace{-0.2cm} \begin{array}{c}
b\\
\varepsilon
\end{array}
 \hspace{-0.2cm}
\right)}
}
)=0$ for $ \varepsilon \in \{1,...,[\frac{k}{2}]\}$.

\underline{Cases   $\varepsilon  \leq  k - 4$.}
We have that every $2$-cycle  $\mu $ in $\Delta_{
\tiny{
\left( \hspace{-0.2cm} \begin{array}{c}
b\\
\varepsilon
\end{array}
 \hspace{-0.2cm}
\right)}
} $
 is homologous    
to a $2$-cycle $\tilde{\mu}$ in $F^{2}(\Delta_b) $ by Lemma \ref{sc1}. 
Since $H_{2}(F^{3}(\Delta_b))=0$
(because, by induction hypothesis,
 $H_{2}(\Delta_b)=0$), we have that 
$[\tilde{\mu}]=0$ in $H_{2}(F^{3}(\Delta_b))=0$.
Since in these cases
 $F^{3}(\Delta_b) 
\subseteq \Delta_{
\tiny{
\left( \hspace{-0.2cm} \begin{array}{c}
b\\
\varepsilon
\end{array}
 \hspace{-0.2cm}
\right)}}$, 
 we may conclude that 
$[\mu ]= [\tilde{\mu}]=0$ in 
 $H_{2}(
\Delta_{
\tiny{
\left( \hspace{-0.2cm} \begin{array}{c}
b\\
\varepsilon
\end{array}
 \hspace{-0.2cm}
\right)}
})$, thus $H_{2}(
\Delta_{
\tiny{
\left( \hspace{-0.2cm} \begin{array}{c}
b\\
\varepsilon
\end{array}
 \hspace{-0.2cm}
\right)}
})=0$.
   
\underline{Cases  $\varepsilon  > k-4$.} 
We have that every $2$-cycle $\mu $
in $ \Delta_{
\tiny{
\left( \hspace{-0.2cm} \begin{array}{c}
b\\
\varepsilon
\end{array}
 \hspace{-0.2cm}
\right)}
}$ is homologous    
to a $2$-cycle $\tilde{\mu}$ in $F^{2}(\Delta_b) $ by Lemma \ref{sc1}. 
Since $H_{2}(F^{3}(\Delta_b))=0$
(because $H_{2}(\Delta_b)=0$), we have that 
$[\tilde{\mu}]=0$ in $H_{2}(F^{3}(\Delta_b))=0$.
But in these cases we have not the inclusion 
 $F^{3}(\Delta_b) 
\subseteq \Delta_{
\tiny{
\left( \hspace{-0.2cm} \begin{array}{c}
b\\
\varepsilon
\end{array}
 \hspace{-0.2cm}
\right)}}$, thus we may not conclude at once. 
Since
$H_{2}(F^{3}(\Delta_b))=0$,
  there exists a $3$-chain $\nu$  in $F^{3}(\Delta_b)$ s.t.
 $\partial \nu= \tilde{\mu}$. We have that
$sp(\nu) =  \cup_{i} F_i$, $F_i$ tetrahedrons in $F^{3}(\Delta_b)$.
Consider a $3$-chain 
$\psi$ in  $\Delta_{\tiny{
\left( \hspace{-0.2cm} \begin{array}{c}
b\\
\varepsilon
\end{array}
 \hspace{-0.2cm}
\right)}
}$ whose support is   
$\cup_{i} \hat{F_i}$, where
$\hat{F_i}$ is a cone $\subseteq \Delta_{\tiny{
\left( \hspace{-0.2cm} \begin{array}{c}
b\\
\varepsilon
\end{array}
 \hspace{-0.2cm}
\right)}
}$ on the border of $F_i$ (there exists by Remark \ref{cone}),
in such way  that $\partial \psi = \tilde{\mu}$;
thus $[\tilde{\mu}] =0 $ in $H_{2}(\Delta_{\tiny{
\left( \hspace{-0.2cm} \begin{array}{c}
b\\
\varepsilon
\end{array}
 \hspace{-0.2cm}
\right)}
})$, thus $[\mu] =0 $ in $H_{2}(\Delta_{\tiny{
\left( \hspace{-0.2cm} \begin{array}{c}
b\\
\varepsilon
\end{array}
 \hspace{-0.2cm}
\right)}
})$.
Then we may conclude that  $H_{2}(
\Delta_{
\tiny{
\left( \hspace{-0.2cm} \begin{array}{c}
b\\
\varepsilon
\end{array}
 \hspace{-0.2cm}
\right)}
})=0$.

\hfill {\it Q.e.d. in Theorem \ref{N2}   }

\section{Proof of Proposition \ref{prodotti}}

Let $X$ and $Y$ be two projective varieties and $L$ a line bundle 
on  $X$ and $M$ a line bundle on $Y$.
Let $\{\sigma_{0},..., \sigma_{k} \}$ be a basis of $H^{0}(X,L)$ and 
let $\{s_{0},..., s_{l} \}$ be a basis of $H^{0}(Y,M)$;
we can suppose $\exists \; \overline{y} \in Y $ s.t. $s_{0}(\overline{y}) 
\neq 0$, $s_{j} (\overline{y})= 0$ for $j \neq 0$;  
let $t_{i,j}$ be the  coordinates corresponding to $\{ 
\sigma_i \otimes s_j\}_{i,j}$
 of the embedding of 
$X \times Y$  by $\pi_X^{\ast}L \otimes \pi_Y^{\ast}M$ 
(where $\pi_{\cdot }$ is the projection on $\cdot $)
and let $t_i$ be the coordinates corresponding to $\{\sigma_{0},...,
 \sigma_{k} \}$ of the embedding of $X$  by $L$.

\begin{remark} \label{eq}
By setting $t_{i,j}=0$ for $j \neq 0$ in an equation of $X \times Y$ and then
taking off the last index (a $0$) of  each variable, we get
an equation of $X$ (to prove this, use $\overline{y}$).  
\end{remark}

\begin{remark} \label{mo}
Let $M$ be a graded module on ${\bf C}[x_1,..., x_n]$ with a minimal set
of generators of degree $s$; then a subset  of elements of degree $s$ of 
$M$ can be extended to a minimal set of 
generators if and only if these elements
are linearly independent on ${\bf C}$.
\end{remark}

{\it Proof of Proposition \ref{prodotti}. }
Suppose $L$ satisfies Property $N_{p-1}$ but not $N_{p}$.
We want to show $\pi_X^{\ast}L \otimes \pi_Y^{\ast}M$ 
does not satisfy  Property $N_p$; we can suppose 
$\pi_X^{\ast}L \otimes \pi_Y^{\ast}M$ 
satisfies  Property $N_{p-1}$. Let $l_m$ and $q_m$
be the ranks of the $m$-module of a minimal free graded 
 resolution respectively of $ G( L) $ 
and of $G(\pi_X^{\ast}L \otimes \pi_Y^{\ast}M) $. 
Let $\{g^{m}_{j}\}_{j=1,..., l_{m}}$, be  a minimal set of generators
of the $m$-module $E_m$ of a minimal resolution,
$... \rightarrow E_m \rightarrow E_{m-1} \rightarrow ....\rightarrow E_0
\rightarrow G(L) \rightarrow 0$,  of $G(L)$.

Since $L$ satisfies Property $N_{p-1}$ but not $N_{p}$, there
exists a  syzygy $S$ of  $(g^{p-1}_{1},...,g^{p-1}_{l_{p-1}})$,
 s.t. $S$ is not generated by linear syzygies of  
$(g^{p-1}_{1},...,g^{p-1}_{l_{p-1}})$.
Add a $0$ to the indices of the variables appearing in $S$  and call
$\tilde{S}$  the so obtained vector  of polynomials; let $\tilde{S}' = 
(\tilde{S}, 0,...,0)$ with $0$ repeated $q_{p-1}- l_{p-1}$ times.
      
Obviously by adding a 
$0$ to the indices of each variable appearing in the equations of $X$, 
we get equations of $X \times Y$ and by adding  a $0$ to the indices of
every variable appearing in the syzygies  of $X$ we get syzygies
of $X \times Y $.

Add a $0$ to the indices of the variables 
appearing in $g^{m}_{j}$  and call
$\tilde{g}^{m}_{j}$ 
the so obtained vector  of  polynomials for  $j=1,..., l_{m}$; set 
$f^{1}_{j}= \tilde{g}^{1}_{j}$ for $j=1,..., l_{1}$ and  
$f^{m}_{j}= (\tilde{g}^{m}_{j},0,...,0)$ 
($0$ repeated $q_{m-1} - l_{m-1}$ times) 
for $j=1,..., l_{m}$ and $2 \leq m \leq p-1$;
$f^{m}_{j}$ for $j=1,..., l_m$ are 
vectors of linear polynomials for  $2 \leq m \leq p-1$ and
 they are quadratic if $m=1$, thus,
by induction on $m$ and  by Remark \ref{mo},
one can extends
this set  to a minimal set of generators
$\{f^{m}_{j}\}_{j=1,..., q_{m}}$, of the $m$-module
of a minimal resolution  of $G(\pi_X^{\ast}L \otimes \pi_Y^{\ast}M) $
 for $m \leq p-1$ (we recall that we supposed  
 $\pi_X^{\ast}L \otimes \pi_Y^{\ast}M $ satisfies Property
 $N_{p-1}$);
we can do it in such way that, when 
we set  $t_{i,j}=0$ for $j \neq 0$, we have that 
$f^{1}_{j}$ is zero for   $j = l_{1} + 1 , ..., q_1$ and 
the $r$-th coordinate of $f^{m}_{j}$ is zero for $r \leq l_{m-1}$ and
 $j = l_{m} + 1 , ..., q_m$
(we can prove this by induction on $m$, by using Remark \ref{eq}
for the case $m=1$: it is sufficient to subtract 
linear combination of  $f^{m}_{j}$  for  $j = 1 , ..., l_m$
 to  $f^{m}_{j}$  for  $j = l_{m} + 1 , ..., q_m$).
 
 Obviously $\tilde{S}'$ is a syzygy of $(f^{p-1}_1,..., f^{p-1}_{q_{p-1}}  )$.

If $\pi_X^{\ast}L \otimes \pi_Y^{\ast}M$ satisfies Property $N_p$
then  $\tilde{S}'$ would be generated by linear syzygies 
 of $(f^{p-1}_1,..., f^{p-1}_{q_{p-1}}  )$.

We state that $\tilde{S}'$ cannot be generated by linear syzygies
 of $(f^{p-1}_1,..., f^{p-1}_{q_{p-1}}  )$.
 In fact, if it were, say $\tilde{S}' = \sum_{\alpha} S_{\alpha}$ 
($S_{\alpha}$
 linear syzygies of $(f^{p-1}_1,..., f^{p-1}_{q_{p-1}}  )$),
 we set  $t_{i,j}=0$ for $j \neq 0$ in each member
of the equality $\tilde{S}' = \sum_{\alpha}
 S_{\alpha}$  and, by taking off the last 
index (a $0$) of  every variable and considering only the first $l_{p-1}$
coordinates of $S$ and $S_{\alpha}$, one would
 obtain that $S$ would be generated by linear 
syzygies of $(g^{p-1}_{1},...,g^{p-1}_{l_{p-1}})$
(observe that by setting  $t_{i,j}=0$ for $j \neq 0$ in $S_{\alpha}$
and taking  the first $l_{p-1}$ coordinates,
we get a syzygy of $(f^{p-1}_{1},...,f^{p-1}_{l_{p-1}})$).
    
But $S$ cannot be generated by linear syzygies  by assumption.      
\hfill {\it Q.e.d.}

\vspace{0.2cm}

By using the program Macaulay \cite{B-S} one can check that  
 ${\cal O}_{{\bf P}^{1} 
\times {\bf P}^{1} \times {\bf P}^{1}}(1,1,1)$ 
  does not satisfy Property $N_4$,  
precisely the resolution, with the notation of Introduction, is:
$$0\rightarrow S(-6) \rightarrow S(-4)^{9} \rightarrow S(-3)^{16}
 \rightarrow
S(-2)^{9} \rightarrow S \rightarrow G \rightarrow 0. $$
From this and from Prop. \ref{prodotti}
we deduce that
 ${\cal O}_{{\bf P}^{1} \times ...\times {\bf P}^{1} }(1,..., 1)$ 
($d$ times) does not satisfy Property $N_4$ for $d \geq 3$. 
By using also Gallego-Purnapranja's Theorem \ref{GP}, we deduce that, 
if  $a_1 , ..., a_d $ are integer numbers
with $a_1 \leq  a_2 \leq ...\leq a_d $ and  $a_{1}=....=a_{k}=1$,
the line bundle ${\cal O}_{{\bf P}^{1} 
\times ...\times {\bf P}^{1} }(a_{1},..., a_{d})$ 
does not satisfy Property $N_4$ if $k \geq 3$
and it does not satisfy Property $N_{2a_{k+1}+2a_{k+2} -2}$ 
if $d- k \geq 2$.
 
With the same argument as in Remark 
in Section 2 of part II of \cite{Green1} we deduce 
 Corollary \ref{N4}.

\vspace{0.7cm}

{\bf Acknowledgements.} I warmly thank G. Ottaviani for suggesting me 
 the problem and for several helpful discussions  and suggestions.
I thank also G. Anzidei and C. Brandigi for some useful discussions.

\end{document}